\documentclass[12pt]{article}
\usepackage{color}
\usepackage{amssymb,latexsym,bm, amsmath,amsthm, amsfonts,multirow, graphicx,mathrsfs}
\usepackage{eufrak,booktabs, ulem}
\usepackage{epsfig}
\usepackage{enumitem}
\graphicspath{{fig/}}

\def\Wb{{\boldsymbol W}}
\def\Xb{{\boldsymbol X}}

\def\Zb{{\boldsymbol Z}}

\def\1{{\bm 1}}
\def\0{{\bm 0}}
\def\th{{\rm th}}

\def\rmspan{{\rm span}}
\def\I{\mathcal{I}}

\def\Sr{{\cal S}_4}

\def\Sr{\mathcal{S}_r}

\renewcommand{\tilde}{\widetilde}
\renewcommand{\hat}{\widehat}

\def\blue{\color{blue}}
\def\red{\color{red}}

\newtheorem{prop}{Proposition}
\newtheorem{thm}{Theorem}

\newtheorem{lemma}{Lemma}

\def\boxit#1{\vbox{\hrule\hbox{\vrule\kern6pt\vbox{\kern6pt#1\kern6pt}\kern6pt\vrule}\hrule}}

\begin{document}

\title{On the asymptotic properties of product-PCA under the high-dimensional setting}

\date{}

\author{Hung Hung$^1$, Chi-Chun Yeh$^1$, and Su-Yun Huang$^2$ \\[0ex]
\small $^1$Institute of Health Data Analytics and Statistics,
National
Taiwan University, Taiwan\\
\small $^2$Institute of Statistical Science, Academia Sinica, Taiwan
}

\maketitle

\begin{abstract}
Principal component analysis (PCA) is a widely used dimension reduction method, but its performance is known to be non-robust to outliers. Recently, product-PCA (PPCA) has been shown to possess the efficiency-loss free ordering-robustness property: (i) in the absence of outliers, PPCA and PCA share the same asymptotic distributions; (ii), in the presence of outliers, PPCA is more ordering-robust than PCA in estimating the leading eigenspace.
PPCA is thus different from the conventional robust PCA methods, and may deserve further investigations. In this article, we study the high-dimensional statistical properties of the PPCA eigenvalues via the techniques of random matrix theory. In particular, we derive the critical value for being distant spiked eigenvalues, the limiting values of the sample spiked eigenvalues, and the limiting spectral distribution of PPCA. Similar to the case of PCA, the explicit forms of the asymptotic properties of PPCA become available under the special case of the simple spiked model. These results enable us to more clearly understand the superiorities of PPCA in comparison with PCA. Numerical studies are conducted to verify our results.

\end{abstract}

\newpage

\section{Introduction}\label{sec.introduction}

Principal component analysis (PCA) is a widely used dimension reduction method. Let $X\in \mathbb{R}^p$ be a random vector with $E(X) = 0$ and $\mathrm{cov}(X) = \Sigma$. The eigenvalue decomposition of $\Sigma$ is given by $\Sigma = \Gamma \Lambda \Gamma^\top$, where
$\Gamma$ is the matrix of eigenvectors and $\Lambda$ is the diagonal matrix of eigenvalues.
Consider the generalized spiked model (GSM), where the eigenvalues satisfy the condition:
\begin{eqnarray}
\lambda_1>\cdots>\lambda_r\gg \lambda_{r+1}\ge \cdots \ge \lambda_p,
\end{eqnarray}
indicating that there exists a rank $r$ such that $\lambda_r$ is significantly larger than $\lambda_{r+1}$. Assuming $r$ is known, a primary goal of PCA is to estimate the leading rank-$r$ eigenspace:
\begin{eqnarray}
\Sr={\rmspan}(\gamma_1,\ldots,\gamma_r),
\end{eqnarray}
and project the data onto $\Sr$ for subsequent analysis. Given the random sample $\{X_i\}_{i=1}^n$, the conventional PCA starts from the sample covariance matrix $\widehat S=\frac{1}{n}\sum_{i=1}^{n}X_iX_i^\top$
to conduct the eigenvalue decomposition $\widehat{S} = \widetilde{\Gamma} \widetilde{\Lambda} \widetilde{\Gamma}^\top$. PCA then estimates $\mathcal{S}_r$ by $\mathrm{span}(\widetilde{\gamma}_1, \ldots, \widetilde{\gamma}_r)$, where the order of $\widetilde{\gamma}_j$ is determined by the order of $\widetilde{\lambda}_j$.

A major drawback of PCA is its non-robustness. Recently, Hung and Huang (2023) develop product-PCA (PPCA) based on the ideas of random-partition and integration. PPCA starts by randomly dividing the data $\{X_i\}_{i=1}^n$ into two subsets to obtain the corresponding covariance matrices $\widehat{S}_1$ and $\widehat{S}_2$. It then calculates the product covariance estimator $\widehat S_{12}=\widehat S_1^{\frac{1}{2}}\widehat S_2^{\frac{1}{2}}$, where $\widehat S_k^{\frac{1}{2}}$ is the positive square root matrix of $\widehat S_k$.
Consider the following singular value decomposition (SVD):
\begin{eqnarray}
\widehat S_{12}=\widehat U\widehat\Lambda \widehat V^\top\quad{\rm with}\quad \widehat U=[\widehat u_1,\ldots,\widehat u_p],\quad
\widehat V=[\widehat v_1,\ldots,\widehat v_p],\quad \widehat \Lambda={\rm diag}(\widehat \lambda_1,\ldots,\widehat \lambda_p). \label{ppca.svd}
\end{eqnarray}
By noting that $\widehat S_{12}$ is an estimator of $\Sigma$, PPCA then estimates $\lambda_j$ by $\widehat \lambda_j$, and estimates $\gamma_j$ by $\widehat\gamma_j = \frac{\widehat u_j+\widehat v_j}{\|\widehat u_j+\widehat v_j\|}$. A basis of $\Sr$ is then estimated by $\widehat\Gamma_r=[\widehat\gamma_1,\ldots,\widehat\gamma_r]$, where the order of $\widehat\gamma_j$ is determined by the order of $\widehat\lambda_j$. PPCA is shown to possess the efficiency-loss free ordering-robustness property: (i) in the absence of outliers, PPCA and PCA share the same asymptotic distribution; (ii), in the presence of outliers, PPCA is more ordering-robust than PCA in estimating the leading eigenspace $\Sr$, indicating that PPCA can perform as well as, or even better than PCA, regardless of the presence of outliers. This characteristic sets PPCA apart from conventional robust PCA methods, which often sacrifice efficiency for robustness gains. This distinction warrants further investigation.

{In Hung and Huang (2023), the authors show that PPCA and PCA share the same asymptotic distribution under either the assumption of finite $p$ with bounded spiked eigenvalues (i.e., $\lambda_1<\infty$) or the assumption of diverging $p$ with diverging spiked eigenvalues (i.e., $\lambda_r\to\infty$). Diverging spiked eigenvalues are easier to be separated from the tail eigenvalues, which makes the estimation problem and the corresponding theoretical investigation easier. However, the asymptotic behavior of PPCA under a more challenging situation, diverging $p$ with bounded spiked eigenvalues, remains unknown. In the situation of diverging $p$ with bounded spiked eigenvalues, it has been shown via the techniques of random matrix theory (RMT) that conventional PCA is biased in estimating
$\lambda_j$'s. In particular, RMT helps quantifying the bias of the spiked eigenvalue estimators $\{\widetilde \lambda_j\}_{j\le r}$, and the limit of the empirical distribution of the non-spiked eigenvalues $\{\widetilde \lambda_j\}_{j> r}$, thereby providing a more clear understanding of PCA's high-dimensional behavior. The aim of this work is to use RMT to study the high-dimensional behavior of PPCA specifically under the scenario of diverging $p$ with bounded spiked eigenvalues.}

The rest of the paper is organized as below. The high-dimensional properties of PPCA are derived in Section~2. Section~3 studies the robustness of PPCA in comparison with PCA from a high-dimensional point of view. Numerical studies are conducted in Sections~4--5. The paper ends with discussions in Section~6. All the proofs are deferred to Appendix.

\section{Asymptotic Properties of PPCA}

The aim of this section is to investigate the asymptotic properties of PPCA using the techniques of RMT. Define the empirical spectral distribution (ESD) of $\Sigma$ to be
\begin{eqnarray}
F^\Sigma(t)=\frac{1}{p}\sum_{j=1}^{p}\mathcal{I}\{\lambda_j\le t\}.
\end{eqnarray}
We aim to study the limiting behavior of the leading eigenvalues $\{\widehat\lambda_j\}_{j\le r}$ and the limiting spectral distribution
(LSD) of $F^{\widehat\Lambda}$ under the following conditions:
\begin{itemize}
\item[(C1)]
$X=\Sigma^{1/2}Z$, where $Z\in \mathbb{R}^p$ has independent elements with zero mean, variance 1, and finite $4^\th$ moment.

\item[(C2)]
$p/n\to c\in(0,\infty)$.

\item[(C3)]
$F^\Sigma$ converges weakly to $H$ as $p\to \infty$, where $H$ has bounded support. Let $H^2$ denote the limiting distribution of $F^{\Sigma^2}$.
\end{itemize}
In the rest of the discussion, $\Sigma$ is said to follow the generalized spiked model (GSM) when $H$ is an arbitrary distribution satisfying (C3), and said to follow the simple spiked model (SSM) when $H$ has positive probability at a single point only. Let $\mathcal{S}(H)$ be the support of $H$, $u_H=\sup \mathcal{S}(H)$ be the upper bound of $\mathcal{S}(H)$, and $\mu_H=\int tdH(t)$ be the mean of $H$. $\mathcal{I}({\cdot})$ is the indicator function.

\subsection{The case of GSM}

To state the main result, we start with introducing the generalized Mar\v{c}enko-Pastur (MP) distribution $F_{c,H}$ indexed by $(c,H)$. Generally, $F_{c,H}$ has no closed form, and its definition is through its Stieltjes transform $m_{c,H}(z)\stackrel{\Delta}{=}\int\frac{1}{t-z}dF_{c,H}(t)$, which is shown to satisfy the {Silverstein equation}:
\begin{eqnarray}
z = -\frac{1}{\underline{m}_{c,H}}+c\int\frac{t}{1+t\underline{m}_{c,H}}dH(t)\quad{\rm with}\quad \underline{m}_{c,H}(z) = cm_{c,H}(z) + \frac{c-1}{z}.\label{s_equation}
\end{eqnarray}
Given $(c,H)$, define the function
\begin{eqnarray}
\psi_{c,H}(\lambda)=\lambda\left\{1+c\int\frac{t}{\lambda-t}dH(t)\right\},\quad \lambda\notin \mathcal{S}(H).\label{psi_lambda}
\end{eqnarray}
The asymptotic properties of PPCA under the GSM are summarized below. For comparison, the results of PCA established in the literature {(Yao, Zheng, and Bai, 2015)} are also provided.

\begin{thm}\label{thm.ppca}
Assume GSM and conditions (C1)--(C3).\\
1. For the case of PPCA, we have the following results.
\begin{enumerate}
\item[(a)]
The ESD of $\widehat \Lambda$ converges to $\lim_{n\to\infty}F^{\widehat \Lambda}=G_{c,H}$, where $G_{c,H}$ is connected with $F_{2c, F_{2c, H^2}}$ via $G_{c,H}(t)=F_{2c, F_{2c, H^2}}(t^2)$.


\item[(b)]
The eigenvalue $\widehat\lambda_j$ can be separated from $\mathcal{S}(G_{c,H})$ if $\lambda_j$ is larger than the threshold $\lambda^*$, which satisfies the following two conditions:
\begin{eqnarray}\label{lambda_star}
\psi_{2c,H^2}(\lambda^{*2}) \ge \sup \mathcal{S}(F_{2c, H^2}) ~~{\rm and}~~ 1 = 2c\int\left(\frac{t}{t-\psi_{2c,H^2}(\lambda^{*2})}\right)^2dF_{2c,H^2}(t).
\end{eqnarray}

\item[(c)]
For any fixed $j$, we have
\begin{eqnarray*}
\widehat\lambda_j&\stackrel{p}{\to}&\Bigg\{\begin{array}{cc}
\psi_j\stackrel{\Delta}{=}\frac{1}{\lambda_j}\psi_{2c,H^2}(\lambda_j^2),&\quad{\rm if}~\lambda_j>\lambda^* , \\
b\stackrel{\Delta}{=}\frac{1}{\lambda^*}\psi_{2c,H^2}(\lambda^{*2}),&
\quad{\rm if}~\lambda_j\le\lambda^* .
\end{array}
\end{eqnarray*}
\end{enumerate}
2. For the case of PCA, we have the following results.
\begin{enumerate}
\item[(a)]
The ESD of $\widetilde \Lambda$ converges to $\lim_{n\to\infty}F^{\widetilde \Lambda}=F_{c,H}$.

\item[(b)]
The eigenvalue $\widetilde\lambda_j$ can be separated from $\mathcal{S}(F_{c,H})$ if $\lambda_j$ is larger than the threshold $\lambda'$, which satisfies the following two conditions:
\begin{eqnarray*}
\lambda' \ge \sup \mathcal{S}(H) ~{\rm and}~ 1 = c\int\left(\frac{t}{t-\lambda'}\right)^2dH(t).
\end{eqnarray*}

\item[(c)]
For any fixed $j$, we have
\begin{eqnarray*}
\widetilde\lambda_j&\stackrel{p}{\to}&\Bigg\{\begin{array}{cc}
\psi_j'\stackrel{\Delta}{=}\psi_{c,H}(\lambda_j),&\quad{\rm if}~\lambda_j>\lambda' , \\
b'\stackrel{\Delta}{=}\psi_{c,H}(\lambda'),&
\quad{\rm if}~\lambda_j\le\lambda' .
\end{array}
\end{eqnarray*}
\end{enumerate}
\end{thm}

One implication of the above theorem is that the threshold for identifiabe eigenvalues can differ between PPCA and PCA. This suggests that a spiked eigenvalue in PCA may not necessarily be a spiked eigenvalue in PPCA. We will delve deeper into this issue in the next subsection.


In high-dimensional PCA, an important issue is the inconsistency of $\widetilde\lambda_j$, as can be seen in Theorem~\ref{thm.ppca}. It demonstrates that the limiting eigenvalues of PCA satisfy $\psi_j'>\lambda_j$ provided that $\lambda_j>\lambda'$ and $c>0$. A similar inconsistency is observed in PPCA, where Theorem~\ref{thm.ppca} indicates that $\widehat\lambda_j$ converges to $\psi_j>\lambda_j$ for a spiked eigenvalue $\lambda_j$ of PPCA. Despite these inconsistencies in both methods, Yata and Aoshima (2010) found that PPCA exhibits smaller bias than PCA when estimating $\lambda_j$. Interestingly, leveraging techniques from RMT allows us to explicitly observe this phenomenon, as elaborated below.

\begin{thm}[]\label{thm.compare.psi}
Assume GSM and conditions (C1)--(C3). We have the following results.
\begin{enumerate}
\item[\rm (i)]
$\psi_j' > \psi_j > \lambda_j$.

\item[\rm(ii)]
$\psi_j'-\psi_j\in(\frac{1}{2}c\mu_H, c\mu_H)$, where $c= \lim_{n,p\to\infty} p/n$.

\item[\rm(iii)]
$|\psi_j-\lambda_j|=O(cu_H^2/\lambda_j)$ and $|\psi_j'-\lambda_j|=O(cu_H)$.

\end{enumerate}
\end{thm}
\noindent The results~(i)--(ii) confirm that PCA always produces larger bias than PPCA in estimating $\lambda_j$, and the difference between their respective limits is bounded within $(\frac{1}{2}c\mu_H, c\mu_H)$. The result (iii) indicates that under the high-dimensional setting with $c>0$, PCA unavoidably encounters bias issues of the order $O(cu_H)$, while PPCA's estimation bias becomes negligible provided that $\lambda_j/u_H\gg cu_H$, reinforcing the advantageous application of PPCA in large $p$ settings.

We close this section by proposing a consistent estimator for $\lambda_j$:
\begin{eqnarray}\label{consistent_est}
\widehat\lambda_j^{\rm c}=-\frac{1}{\underline{\widehat{s}}_j(\widehat\lambda_j^2)\widehat\lambda_j},
\end{eqnarray}
where $\underline{\widehat{s}}_j(z)=2c\widehat s_j(z)+(2c-1)z^{-1}$ and $\widehat s_j(z)=\frac{1}{p-j}\sum_{\ell > j}(\widehat\lambda_\ell^2-z)^{-1}$. The consistency of $\widehat\lambda_j^{\rm c}$ is established by the following theorem.

\begin{thm}\label{thm.ppca.consistency}
Assume GSM and conditions~(C1)--(C3). Then, for any $\lambda_j>\lambda^*$,  where $\lambda^*$ is a constant satisfying conditions given by~(\ref{lambda_star}), the estimator $\widehat\lambda_j^{\rm c}$ in~(\ref{consistent_est}) converges in probability to $\lambda_j$ as $n\to\infty$.
\end{thm}

\noindent Note also that a consistent estimator $\widetilde\lambda_j^{\rm c}=-\frac{1}{\underline{\widetilde m}_j(\widetilde\lambda_j)}$ based on PCA is proposed by Bai and Ding (2012), where $\underline{\widetilde{m}}_j(z)=c\widetilde{m}_j(z)+(c-1)z^{-1}$ and $\widetilde m_j(z)=\frac{1}{p-j}\sum_{\ell > j}(\widetilde\lambda_\ell-z)^{-1}$. We will further compare the performances of $\widehat\lambda_j^{\rm c}$ and $\widetilde\lambda_j^{\rm c}$ via numerical studies in Section~\ref{sec.sim}.

\subsection{The case of SSM}\label{sec.ssm}

An issue for RMT is that the knowledge of $F_{c,H}$ is only characterized via $m_{c,H}$, and the explicit form of $F_{c,H}$ is rarely known for general $H$. One exception is the case of SSM, where $H(t)=\mathcal{I}\{t\ge\sigma^2\}$ for some $\sigma^2>0$. In this situation, (\ref{psi_lambda}) becomes $\psi_{c,\sigma^2}(\lambda)=\lambda(1+\frac{c\sigma^2}{\lambda-\sigma^2})$, and $F_{c,H}$ reduces to the MP distribution (denoted by $F_{c,\sigma^2}$) whose pdf has a closed form. Notably, the same situation applies to PPCA, where the LSD of $\widehat\Lambda$, denoted by $G_{c,\sigma^2}$, (where $G_{c,\sigma^2}(t)=F_{2c, F_{2c, \sigma^4}}(t^2)$ as defined in Theorem~\ref{thm.ppca}), also has a closed form as stated below.

\begin{thm}\label{thm.ppca.ssm}
Assume SSM and conditions (C1)--(C3).\\
1. For PPCA, we have the following results.
\begin{enumerate}
\item[(a)]
The ESD of $\widehat \Lambda$ converges to $\lim_{n\to\infty}F^{\widehat \Lambda}=G_{c,\sigma^2}$, which has the pdf:
\begin{eqnarray*}
g_{c,\sigma^2}(t)
&=&\frac{1}{\kappa(t)\pi c\sigma^2}\sqrt{(t^2-\alpha)(\beta-t^2)}\cdot\mathcal{I}\{\alpha\le t^2\le \beta\}
\end{eqnarray*}
with an additional point mass at $0$ with probability $1-(2c)^{-1}$  if $c>1/2$, where
\begin{eqnarray*}
\alpha=\frac{(2+10c-c^2)-\sqrt{c(c+4)^3}}{2}\sigma^4,\quad \beta=\frac{(2+10c-c^2)+\sqrt{c(c+4)^3}}{2}\sigma^4,
\end{eqnarray*}
and $\kappa(t)>0$ is defined as
\begin{eqnarray*}
\kappa(t)&=&\left\{t
\sqrt{(t^2-\alpha)(\beta-t^2)}+\frac{9(c+1)\sigma^2t^2+(2c-1)^3\sigma^6}{3\sqrt{3}}\right\}^{\frac{2}{3}}\\
&&+\left\{t
\sqrt{(t^2-\alpha)(\beta-t^2)}-\frac{9(c+1)\sigma^2t^2+(2c-1)^3\sigma^6}{3\sqrt{3}}\right\}^{\frac{2}{3}}+\frac{3t^2+(2c-1)^2\sigma^4}{3}.
\end{eqnarray*}
Moreover, $\mathcal{S}(G_{c,\sigma^2})=\left[a, b\right]$, where
\begin{eqnarray*}
a&=&\sigma^2\left(1+c-\sqrt{c^2+4c}\right)^{\frac{1}{2}}\left(1-\frac{\sqrt{c^2+4c}+c}{2}\right)\cdot \mathcal{I}\{c<1/2\},\\
b&=&\sigma^2\left(1+c+\sqrt{c^2+4c}\right)^{\frac{1}{2}}\left(1+\frac{\sqrt{c^2+4c}-c}{2}\right).
\end{eqnarray*}

\item[(b)]
The eigenvalue $\widehat\lambda_j$ can be separated from $\mathcal{S}(G_{c,\sigma^2})$ if $\lambda_j$ is larger than the threshold $\lambda^*=\sigma^2\left(1+c+\sqrt{c^2+4c}\right)^{1/2}$.

\item[(c)]
For any fixed $j$, we have
\begin{eqnarray*}
\widehat\lambda_j&\stackrel{p}{\to}&\Bigg\{\begin{array}{cc}
\lambda_j\left(1+\frac{2c\sigma^4}{\lambda_j^2-\sigma^4}\right),&\quad{\rm if}~\lambda_j>\lambda^*,  \\
b,&
\quad{\rm if}~\lambda_j\le\lambda^*.
\end{array}
\end{eqnarray*}
\end{enumerate}
2. For PCA, we have the following results.
\begin{enumerate}
\item[(a)]
The ESD of $\widetilde \Lambda$ converges to $\lim_{n\to\infty}F^{\widetilde \Lambda}=F_{c, \sigma^2}$, which has the pdf:
\begin{eqnarray*}
f_{c,\sigma^2}(t)=\frac{1}{2t\pi c\sigma^2} \sqrt{(t-a')(b'-t)}\cdot\I\{a'\le t\le b'\}
\end{eqnarray*}
with an additional point mass at $t=0$ with probability $1-c^{-1}$ if $c>1$, where $a'=\sigma^2(1-\sqrt{c})^2$, $b'=\sigma^2(1+\sqrt{c})^2$. Moreover, $\mathcal{S}(F_{c,\sigma^2})=[a',b']$.

\item[(b)]
The eigenvalue $\widetilde\lambda_j$ can be separated from $\mathcal{S}(F_{c,\sigma^2})$ if $\lambda_j$ is larger than the threshold $\lambda'=\sigma^2(1+\sqrt{c})$.

\item[(c)]
For any fixed $j$, we have
\begin{eqnarray*}
\widetilde\lambda_j&\stackrel{p}{\to}&\Bigg\{\begin{array}{cc}
\lambda_j(1+\frac{c\sigma^2}{\lambda_j-\sigma^2}),&\quad{\rm if}~\lambda_j>\lambda' ,\\
b',&
\quad{\rm if}~\lambda_j\le \lambda'.
\end{array}
\end{eqnarray*}
\end{enumerate}
\end{thm}
\noindent
\noindent {Figure~\ref{fig.ppca_pdf} reports the histograms of the empirical eigenvalues of PPCA obtained from a random sample generated from $N(0,I_p)$ with $n=2000$ and $p\in\{800, 4000\}$, which gives $c\in\{0.4, 2\}$ (the case of PCA is also reported for comparison). It can be seen that the histogram perfectly matches the derived form of $g_{c,1}$ as expected.}

From the explicit expressions in Theorem~\ref{thm.ppca.ssm}, we have the following observations regarding the performance differences between PPCA and PCA.

\begin{figure}[!b]
\centering
\includegraphics[height=7cm]{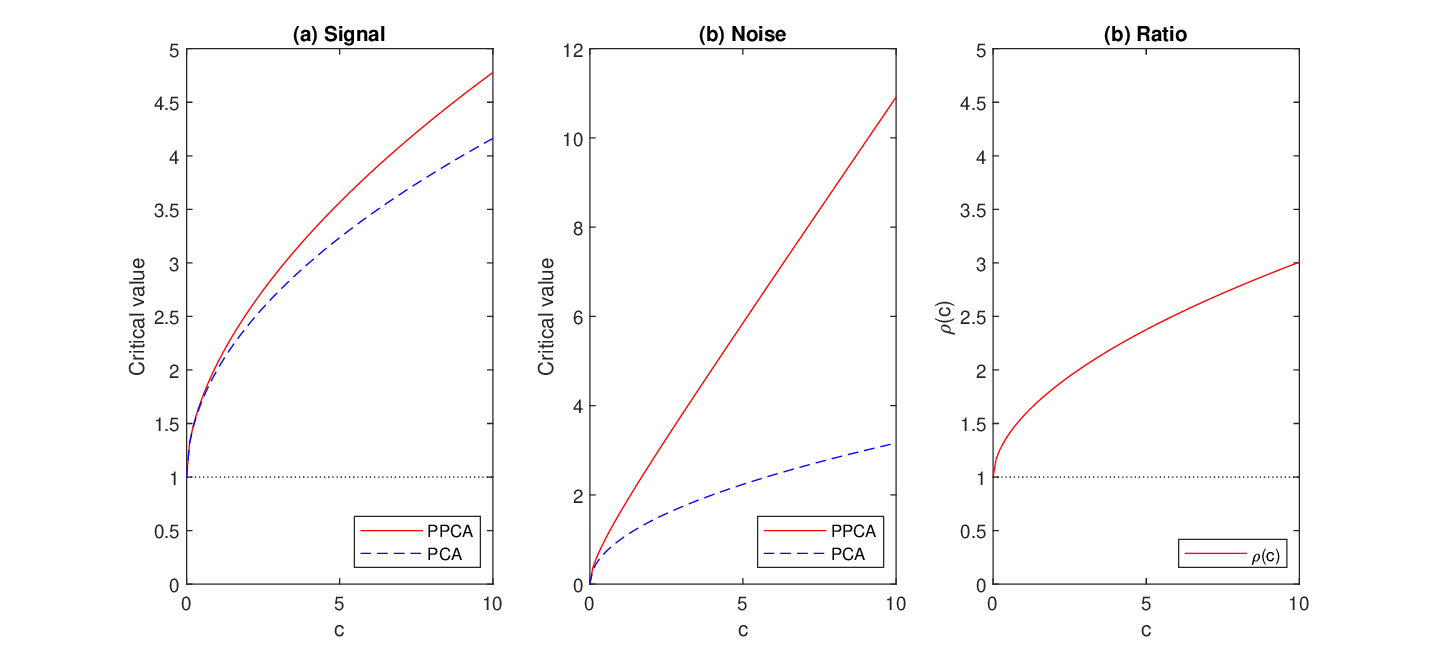}
\caption{Comparison of PPCA and PCA under SSM with $\sigma^2=1$ at $c\in [0,10]$. (a) The thresholds $\lambda^*$ of PPCA and $\lambda'$ of PCA for spiked signal eigenvalues. (b) The critical values $\frac{c+\sqrt{c^2+4c}}{2}$ of PPCA and $\sqrt{c}$ of PCA for spiked noise eigenvalues. (c) The values of the ratio $\rho(c)$. The dashed line represents the value of 1.}
\label{fig.rho_c}
\end{figure}

\begin{itemize}
\item
In PCA, the lower bound $a'$ of $\mathcal{S}(F_{c,\sigma^2})$ increases with $|\sqrt{c}-1|$, indicating that
PCA tends to produce {larger non-zero} empirical tail eigenvalues in both cases $p\gg n$ and $p\ll n$. However, the situation is different for empirical eigenvalues produced by PPCA. In PPCA, the lower bound $a$ of $\mathcal{S}(G_{c,\sigma^2})$ can be strictly positive only when $c<1/2$, and the lower bound $a=0$ when $c\ge 1/2$. For $a=0$,  it indicates that PPCA tends to produce {smaller} empirical tail eigenvalues, which helps in separating empirical tail eigenvalues from spiked eigenvalues.

\begin{figure}[h]
\centering
\includegraphics[height=6cm]{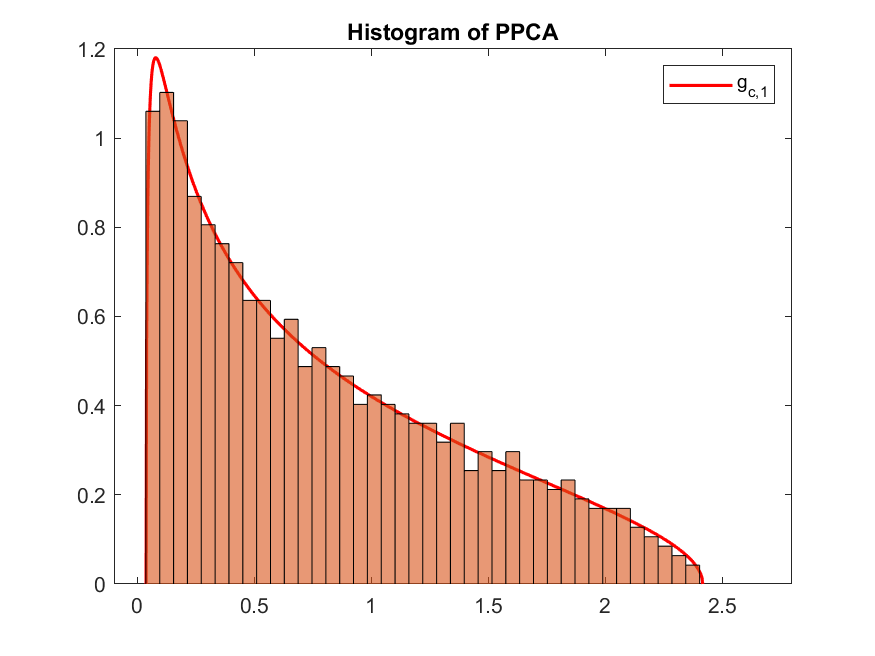}
\includegraphics[height=6cm]{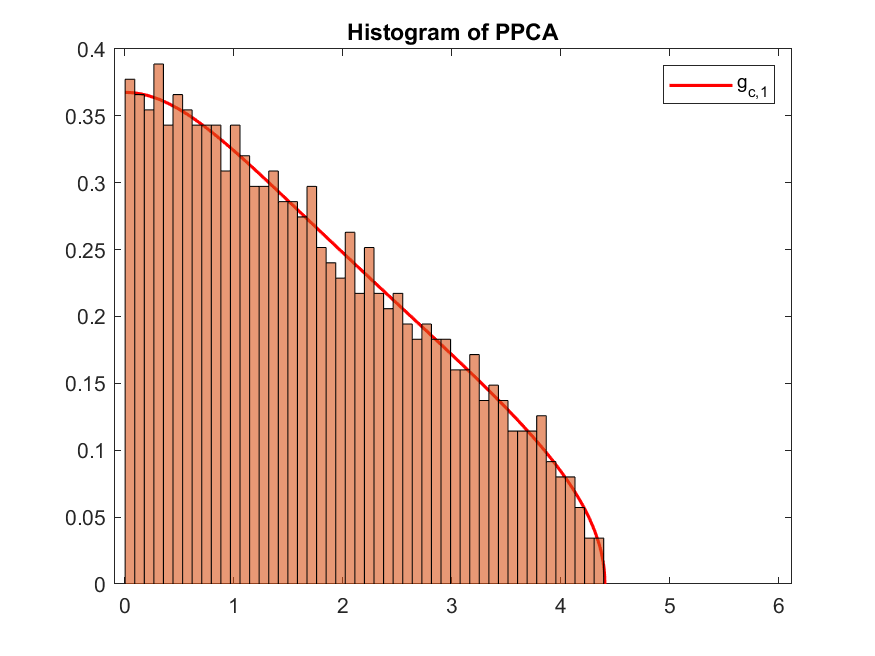}
\includegraphics[height=6cm]{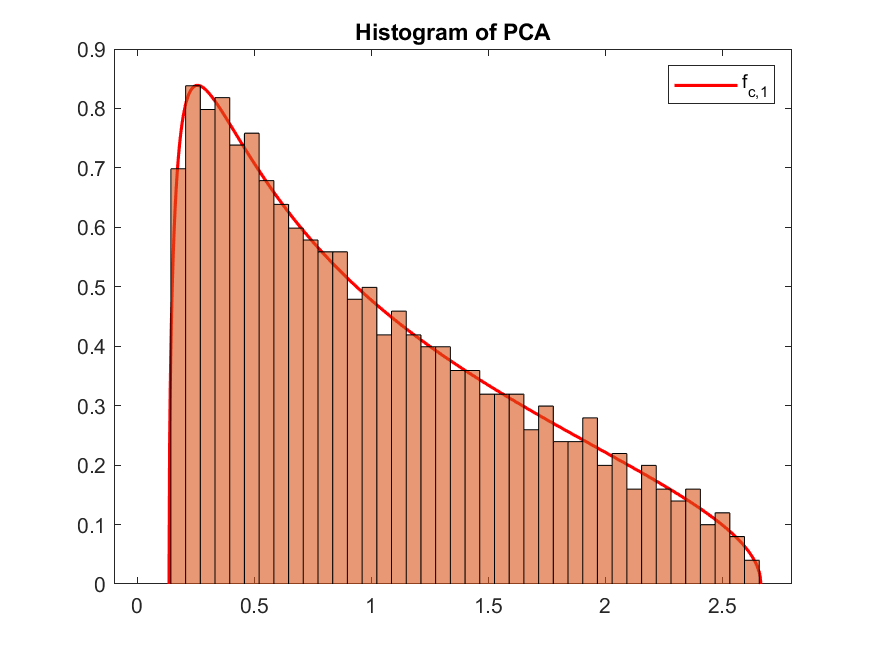}
\includegraphics[height=6cm]{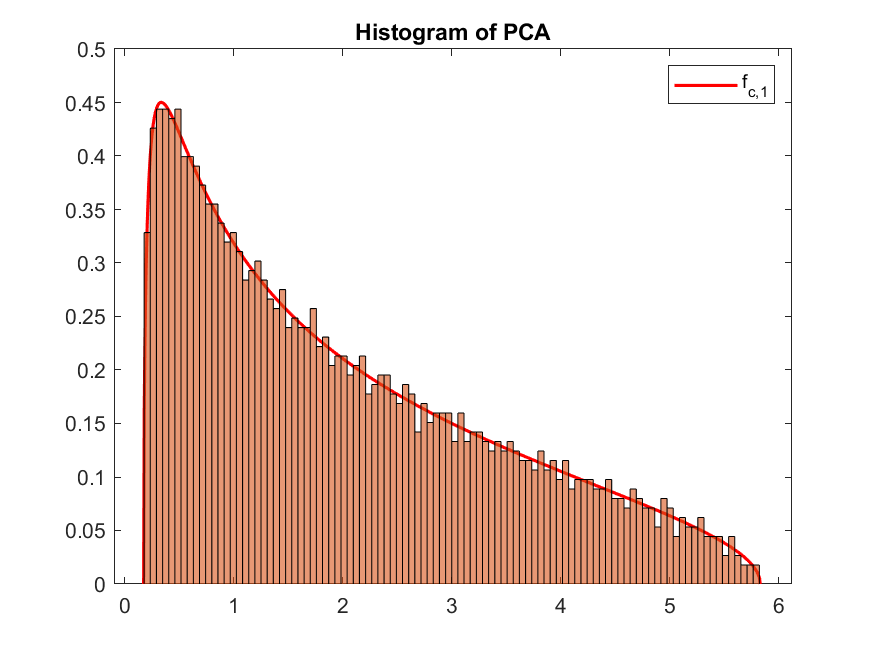}
\caption{The histograms of $\{\widehat\lambda_j:1\le j \le \min\{n/2,p\}\}$ from PPCA and $\{\widetilde\lambda_j:1\le j \le \min\{n,p\}\}$ from PCA under one realization of $X_1,\ldots,X_n\sim N(0,I_p)$ with $n=2000$ and $p\in\{800,4000\}$ that correspond to $c=0.4$ (the left panel) and $c=2$ (the right panel). The histogram is normalized such that the area matches $g_{c,1}$ for PPCA and $f_{c,1}$ for PCA, and the plots of $g_{c,1}$ and $f_{c,1}$ are also reported for comparisons.}

\label{fig.ppca_pdf}
\end{figure}

\item
Under SSM, Theorems~\ref{thm.ppca} and~\ref{thm.compare.psi} lead to
\begin{eqnarray*}
\lim_{n\to\infty}|\widehat\lambda_j-\lambda_j|=O(c\sigma^4/\lambda_j)\quad{\rm and}\quad \lim_{n\to\infty}|\widetilde\lambda_j-\lambda_j|=O(c\sigma^2).
\end{eqnarray*}
Thus, the bias of the PPCA eigenvalue estimate is ignorable when $\lambda_j/\sigma^2 \gg c\sigma^2$, while PCA always suffers a bias issue of order $O(c\sigma^2)$. Moreover, the improvement of PPCA can be explicitly seen from the expression of the limiting value:
\begin{eqnarray*}
\psi_j = \lambda_j\left(1+\frac{c\sigma^2}{\lambda_j-\sigma^2}\cdot\frac{2\sigma^2}{\lambda_j+\sigma^2}\right).
\end{eqnarray*}
Comparing this to the limiting value for PCA
\[\psi_j'=\lambda_j\left( 1+\frac{c\sigma^2}{\lambda_j-\sigma^2}\right),\]
PPCA reduces the bias of eigenvalue estimates. Specifically, PPCA's bias is smaller by a factor of $\frac{2\sigma^2}{\lambda_j+\sigma^2}<1$ compared to PCA.

\item
Comparing the critical values $\lambda^*=(1+c+\sqrt{c^2+4c})^{1/2}$ of PPCA with $\lambda'=1+\sqrt{c}$ of PCA, it is evident that $\lambda^*>\lambda'$ for any $c>0$. Therefore, PPCA may fail to identify eigenvalues in the interval $[\lambda',\lambda^*]$, whereas PCA is capable of identifying eigenvalues within this interval. This highlights a limitation of PPCA, although the difference $\lambda^*-\lambda'$ is small for moderate $c$ {(see Figure~\ref{fig.rho_c}~(a) for the curves of $\lambda^*$ and $\lambda'$ at different $c$ values)}. Moreover, the failure to identify these eigenvalues comes with a significant benefit for PPCA, as it enables the method to avoid incorrectly treating outliers as signal eigenvalues. This issue will be further discussed in Section~\ref{sec.robustness}.
\end{itemize}

\section{Ordering-Robustness of PPCA}\label{sec.robustness}

The aim of this section is to investigate the robustness of PPCA from a high-dimensional point of view and to compare it with PCA. In this situation, the influence of outliers is twofold:
\begin{itemize}
\item[(i)]
Outliers can become visible and be wrongly included as spiked eigenvalues (referred to as spiked noise eigenvalues), which can lead to an overestimation of the target rank $r$ when using any rank selection method (e.g., AIC or BIC). In this situation, even if the ordering of sample eigenvalues remains unaffected, the estimation of $\Sr$ is at risk of erroneously including a noise eigenvector associated with a spiked noise eigenvalue. This inclusion results in the overestimation of $\Sr$, producing a larger subspace.

\item[(ii)]
If an outlier is significant enough to not only become a spiked noise eigenvalue but also to affect the ordering of sample eigenvalues by surpassing certain signal eigenvalues $\{\lambda_j\}_{j\le r}$, the estimation of $\Sr$ becomes biased even if $r$ is correctly specified. In this situation, the leading $r$-dimensional eigenspace will include the eigenvector associated with the spiked noise eigenvalue.
\end{itemize}
\noindent In Hung and Huang (2023), the authors only discuss issue (ii) by comparing the influence functions of PPCA and PCA at the population level. They consider scenarios with arbitrary but finite $p$, and under a simple outlier model featuring a single outlier.
Under the finite $p$ assumption, every eigenvalue is a spiked eigenvalue. In this article, we aim to study issues~(i) and~(ii) under a more general outlier model. In the rest of discussion, a ``spiked eigenvector'' is an eigenvector whose corresponding eigenvalue is a spiked eigenvalue.

Consider the following SSM, assuming $\sigma^2=1$ without loss of generality:
\begin{eqnarray}
\Sigma=\lambda_1\gamma_1\gamma_1^\top + (I-\gamma_1\gamma_1^\top),
\end{eqnarray}
where $\lambda_1$ is a single signal eigenvalue satisfying $\lambda_1>(1+c+\sqrt{c^2+4c})^{1/2}$. From Theorem~\ref{thm.ppca.ssm}, $\lambda_1$ qualifies as a spiked eigenvalue for both PCA and PPCA, leading to $\Sr=\rmspan(\gamma_1)$ with $r=1$. Suppose that in the presence of outliers $\Sigma$ is perturbed to
\begin{eqnarray}
\Sigma_\varepsilon=(1-K\epsilon)\Sigma+\epsilon\sum_{k=1}^{K}\eta_k\nu_k\nu_k^\top,\label{Sigma.perturbed}
\end{eqnarray}
where $\eta_k>0$ is the effect size and $\nu_k$ represents the unit-length direction of the $k^\th$ perturbation, $k=1,\ldots,K$. The parameter $\epsilon>0$ controls the proportion of contamination.
Assuming $\nu_k$'s are random directions in a high-dimensional space, it is reasonable to expect that, for any $k\ne\ell$, $\gamma_1^\top\nu_k= 0$ and $\nu_k^\top\nu_\ell=0$ hold approximately.

{\bf PCA.} Direct calculation gives the population PCA with the following spectral decomposition:
\begin{eqnarray}\label{ssm_pca}
\Sigma^{\rm (pca)}_\epsilon &=&(1-K\epsilon)\lambda_1\gamma_1\gamma_1^\top+\sum_{k=1}^{K}(1-K\epsilon+\epsilon\eta_k)\nu_k\nu_k^\top+(1-K\epsilon)Q,
\end{eqnarray}
where $Q=I-\gamma_1\gamma_1^\top-\sum_{k=1}^{K}\nu_k\nu_k^\top$. {Note that $\Sigma^{\rm (pca)}_\epsilon$ follows the SSM with $\sigma^2 = (1-K\epsilon)$. It then implies from Part-2 of Theorem~\ref{thm.ppca.ssm} that, if the effect size $\eta_k$ is large enough to satisfy $1-K\varepsilon+\varepsilon\eta_k > \lambda'$, or equivalently,
\begin{equation}
\eta_k> \frac{1-K\epsilon}{\epsilon}\sqrt{c}, \label{ex.PCA_eta}
\end{equation}
then $(1-K\epsilon+\epsilon\eta_k)$ becomes a (distant) spiked eigenvalue of $\Sigma^{\rm (pca)}_\epsilon $ with the corresponding eigenvector $\nu_k$}. In this situation, PCA suffers the risk of wrongly including $\nu_k$ in the estimation of $\Sr$, and the target rank of PCA is no longer $r$ but becomes
\begin{eqnarray}
r_{\rm pca}:=r + \sum_{k=1}^K\mathcal{I}\left\{\eta_k>\frac{1-K\varepsilon}{\varepsilon}\sqrt{c}\right\}.\label{r.pca}
\end{eqnarray}
Furthermore, if the effect size $\eta_k$ is large enough to satisfy $1-K\varepsilon+\varepsilon\eta_k > (1-K\varepsilon)\lambda_1$, or equivalently,
\begin{eqnarray}
\eta_k>\frac{1-K\varepsilon}{\varepsilon}(\lambda_1-1),
\label{ex.PCA_a}
\end{eqnarray}
then the ordering of $\nu_k$ will precede the ranking of $\gamma_1$, and the estimation of $\Sr$ based on PCA
is biased.

{\bf PPCA.} Recall that $\Xb$ is randomly partitioned into $\{\Xb_1,\Xb_2\}$ with equal size. Assuming that $\epsilon$ is small enough such that the perturbation $\eta_k\nu_k\nu_k^\top$ can appear in only one of $\{\Xb_1,\Xb_2\}$, this leads to the sub-sample covariance matrices:
\begin{eqnarray*}
&&\Xb_1\to \Sigma_{1,\varepsilon}=(1-2K_1\varepsilon) \Sigma +2\varepsilon\sum_{k\in \mathcal{K}_1}\eta_k\nu_k\nu_k^\top,\\
&&\Xb_2\to \Sigma_{2,\varepsilon}=(1-2K_2\varepsilon) \Sigma + 2\varepsilon\sum_{k\in \mathcal{K}_2}\eta_k\nu_k\nu_k^\top,
\end{eqnarray*}
where $\{\mathcal{K}_j\}_{j=1}^2$ is a partition of $\{1,2,\ldots,K\}$ corresponding to the partition  $\{\Xb_1,\Xb_2\}$ with $K_j=|\mathcal{K}_j|$, and $K_1+K_2=K$. That is, the indices of $\nu_k$'s are randomly allocated to either $\mathcal{K}_1$ or $\mathcal{K}_2$, corresponding to $\Xb_1$ or $\Xb_2$. The factor of 2 adjusts for the proper influence scale, reflecting that the sample size of $\Xb_1$ and $\Xb_2$ is only $n/2$ each. Then, the covariance matrix of PPCA, $\Sigma^{\rm (ppca)}_\varepsilon := \Sigma_{1,\varepsilon}^{\frac{1}{2}}\Sigma_{2,\varepsilon}^{\frac{1}{2}}$, has the following spectral decomposition:
\begin{eqnarray}
\Sigma^{\rm (ppca)}_\varepsilon
&=&\{(1-2\varepsilon K_1)(1-2\varepsilon K_2)\}^{\frac{1}{2}}\lambda_1\gamma_1\gamma_1^\top\nonumber \\
&&+\sum_{k\in \mathcal{K}_1}\{(1-2\varepsilon K_2)(1-2\varepsilon K_1+2\varepsilon\eta_k)\}^{\frac{1}{2}}\nu_k\nu_k^\top\nonumber \\
&&+\sum_{k\in \mathcal{K}_2}\{(1-2\varepsilon K_1)(1-2\varepsilon K_2+2\varepsilon\eta_k)\}^{\frac{1}{2}}\nu_k\nu_k^\top\nonumber \\
&& + \{(1-2\varepsilon K_1)(1-2\varepsilon K_2)\}^{\frac{1}{2}}Q, \label{ssm_ppca}
\end{eqnarray}
which follows the SSM with $\sigma^2=\{(1-2\varepsilon K_1)(1-2\varepsilon K_2)\}^{\frac{1}{2}}$. Consider a case where $k\in \mathcal{K}_1$ (the case of $\mathcal{K}_2$ is similar). {It then implies from Part-1 of Theorem~\ref{thm.ppca.ssm} that, if the size $\eta_k$ is large enough to satisfy $\{(1-2\varepsilon K_2)(1-2\varepsilon K_1+2\varepsilon\eta_k)\}^{\frac{1}{2}} > \lambda^*$, or equivalently,
\begin{equation}
\quad \eta_k>\frac{1-2K_1\varepsilon}{\varepsilon}\cdot\frac{c+\sqrt{c^2+4c}}{2}, \label{ex.PPCA_eta}
\end{equation}
then $\{(1-2\varepsilon K_2)(1-2\varepsilon K_1+2\varepsilon\eta_k)\}^{\frac{1}{2}}$ becomes a (distant) spiked eigenvalue of $\Sigma^{\rm (ppca)}_\varepsilon$ with the corresponding eigenvector $\nu_k$}. In this situation, PPCA suffers the risk of wrongly including $\nu_k$ in the estimation of $\Sr$, and the target rank of PPCA becomes
\begin{eqnarray}
r_{\rm ppca}:= r + \sum_{\ell=1}^2\sum_{k\in \mathcal{K}_\ell}\mathcal{I}\left\{\eta_k>\frac{1-2K_\ell\varepsilon}{\varepsilon}\cdot\frac{c+\sqrt{c^2+4c}}{2}\right\}.\label{r.ppca}
\end{eqnarray}
If $\eta_k$ further satisfies $\{(1-2\varepsilon K_2)(1-2\varepsilon K_1+2\varepsilon\eta_k)\}^{\frac{1}{2}} > \{(1-2\varepsilon K_1)(1-2\varepsilon K_2)\}^{\frac{1}{2}}\lambda_1$, or equivalently,
\begin{eqnarray}
\eta_k>\frac{1-2K_1\varepsilon}{2\varepsilon}(\lambda_1^2-1), \label{ex.PPCA_a}
\end{eqnarray}
then the ordering of $\nu_k$ will precede the ordering of $\gamma_1$, and the estimation of $\Sr$ based on PPCA
is biased.

We have the following observations regarding the robustness comparison of PPCA and PCA:
\begin{itemize}
\item
Comparing conditions on $\eta_k$ for spiked noise eigenvalue, (\ref{ex.PCA_eta}) for PCA and (\ref{ex.PPCA_eta}) for PPCA, it requires a larger size of $\eta_k$ for $\nu_k$ to become a spiked eigenvector in the case of PPCA than in PCA under the condition $\frac{1-2K_1\varepsilon}{\varepsilon}\cdot\frac{c+\sqrt{c^2+4c}}{2}>\frac{1-K\varepsilon}{\varepsilon}\sqrt{c}$ {(see also Figure~\ref{fig.rho_c}~(b) for two critical values $\sqrt{c}$ and $\frac{c+\sqrt{c^2+4c}}{2}$ at different $c$ values)}, or equivalently,
\begin{eqnarray}
\frac{1-2K_1\varepsilon}{1-K\varepsilon}\cdot\frac{\sqrt{c}+\sqrt{c+4}}{2} > 1.\label{eta_win}
\end{eqnarray}
Note that $\frac{\sqrt{c}+\sqrt{c+4}}{2}>1$ for any $c>0$ and (\ref{eta_win}) holds for small $\varepsilon$, and in this situation, the extra-inclusion dimensions $\delta_{\rm pca}$ and $\delta_{\rm ppca}$ satisfy the inequality $\delta_{\rm pca}\ge\delta_{\rm ppca}$. This indicates that PPCA tends to incur less risk than PCA in incorrectly including $\nu_k$'s in the estimation of $\Sr$. This analysis shows that PPCA exhibits greater resistance to perturbation $\nu_k$ becoming a spiked eigenvector, which is advantageous since $\nu_k$ is typically regarded as contamination or an outlier direction.

\item
Comparing conditions for $\lambda_1$ being leading spiked eigenvalue, (\ref{ex.PCA_a}) for PCA and (\ref{ex.PPCA_a}) for PPCA, it requires a larger size of $\eta_k$ to precede the ranking of $\lambda_1$ in PPCA than in PCA under the condition $\frac{1-2K_1\varepsilon}{2\varepsilon}(\lambda_1^2-1)>\frac{1-K\varepsilon}{\varepsilon}(\lambda_1-1)$, or equivalently,
\begin{eqnarray}
\lambda_1>\frac{1-2K_2\varepsilon}{1-2K_1\varepsilon}. \label{a_win}
\end{eqnarray}
Note that (\ref{a_win}) holds for small $\varepsilon$ since $\lambda_1>1$. This indicates that, in the presence of perturbation $\nu_k$, PPCA demonstrates better capability in maintaining the correct ordering between $\nu_k$ and $\gamma_1$. Ensuring the correct ordering of eigenvectors is crucial for unbiased estimation of $\Sr$.

\item
{There exists a trade-off between the identification range for signal and noise eigenvalues. In particular, PPCA fails to identify signal eigenvalues in $[\lambda',\lambda^*]$ (see Figure~\ref{fig.rho_c}(a)) but gains a larger space to prevent noise eigenvalues from being spiked eigenvalues than PCA (see Figure~\ref{fig.rho_c}(b)).} To investigate the gain and loss of PPCA, define
\begin{eqnarray}
\rho(c) = \frac{(\sqrt{c}+\sqrt{c+4})/2}{\lambda^*/\lambda'}
=\frac{1+\sqrt{c}}{\sqrt{2}}\left(\frac{2+c+\sqrt{c^2+4c}}{1+c+\sqrt{c^2+4c}}\right)^{\frac{1}{2}},
\end{eqnarray}
where $\lambda^*/\lambda'$ is the ratio of threshold values of PPCA and PCA for spiked signal eigenvalues, and $(\sqrt{c}+\sqrt{c+4})/2$ controls the ratio of threshold values for spiked noise eigenvalues in~(\ref{eta_win}). A value of $\rho(c)>1$ then indicates that PPCA gains more space to tolerate noise eigenvalues than the loss of space to identify signal eigenvalues. We have the following results (see also Figure~\ref{fig.rho_c}(c) for the curve of $\rho(c)$ at different $c$ values).

\begin{prop}\label{prop.rho}
The function $\rho(c)$ increases with $c$ and satisfies $\rho(c)\ge 1$ for $c\ge 0$.
\end{prop}

\noindent This result supports the usage of PPCA, especially for the case of moderate to strong spiked signal eigenvalues. {In this situation, both PPCA and PCA are capable of identifying $\Sr$, but PPCA tends to produce less noise spiked eigenvalues than PCA does when outliers are present. This fact also echoes the efficiency-loss free ordering-robustness property of PPCA established in Hung and Huang (2023) under the high-dimensional setting with diverging signal eigenvalues.}

\item
When $K_1=K/2$, (\ref{eta_win})--(\ref{a_win}) automatically hold no matter the sizes of $\varepsilon$ and $c$. This suggests implementing PPCA using a random-partition with $n_1=n_2=n/2$, in this situation the case of $K_1=K/2$ occurs with the highest probability ${K\choose K/2}2^{-K}$.

\item
The most severe situation for PPCA corresponds to the case of $K_1=K$, under which (\ref{eta_win})--(\ref{a_win}) become
\begin{eqnarray}
\frac{1-2K\varepsilon}{1-K\varepsilon}\cdot\frac{\sqrt{c}+\sqrt{c+4}}{2} > 1 \quad{\rm and}\quad \lambda_1>\frac{1}{1-2K\varepsilon}.\label{ppca_K}
\end{eqnarray}
Despite in this situation, (\ref{ppca_K}) can still be satisfied when $\varepsilon$ is small. Note that a small $\varepsilon$ is a reasonable assumption due to the nature of outliers. Note also that under the random-partition with $n_1=n_2=n/2$, the case of $K_1=K$ occurs with the lowest probability $2^{-K}$.
\end{itemize}

The above discussion not only supports the eigenvalue-ordering robustness of PPCA but also indicates that an equal partition, i.e., $n_1=n_2=n/2$, is the optimal size of random-partition to implement PPCA. This conclusion is also found in Yata and Aoshima (2010), where the authors demonstrate that the choice $n_1=n_2=n/2$ achieves the minimum asymptotic variance (in the scenario of absence of outliers). From a robustness standpoint, our discussion further underscores the optimality of the choice $n_1=n_2=n/2$ in situations where outliers are present.

\section{Simulation Studies}\label{sec.sim}

\subsection{Asymptotic behavior of PPCA}

For each simulation replicate, {we generate $\Gamma$ as the orthogonal basis of a random matrix with elements being generated from the standard normal distribution}, and set the diagonal elements of $\Lambda$ to be
\begin{eqnarray}
{\rm diag}(\Lambda)=(10\lambda^*, 5\lambda^*, 1,\ldots,1),
\end{eqnarray}
where $\lambda^*=(1+c+\sqrt{c^2+4c})^{1/2}$ with $c=p/n$. This implies that the first two eigenvalues are identifiable by both PPCA and PCA. Hence, the target rank of $\Sigma$ is $r=2$, which gives $\Sr=\rmspan([\gamma_1,\gamma_2])$. Given $(\Gamma,\Lambda)$, the simulation data $\{X_i\}_{i=1}^n$ is then generated from $X_i=(\frac{\nu-2}{\nu}\Sigma)^{1/2}Z_i$, where the elements of the random vector $Z_i\in \mathbb{R}^p$ are generated from the t-distribution $t_\nu$ with degrees-of-freedom $\nu$. We then conduct PPCA to report the simulation results (with 200 replicates) under $\nu=30$, $n=500$, and $p\in\{200,1000\}$, i.e., $c\in\{0.4,2\}$. Table~\ref{table.1} reports the means of $(\widehat\lambda_1, \widehat\lambda^{\rm c}_1,\widehat\lambda_2, \widehat\lambda^{\rm c}_2, \widehat\lambda_3, \widehat\lambda_p)$ from PPCA in comparison with $(\psi_1,\lambda_1,\psi_2,\lambda_2,b,a)$.

\begin{table}[h]\label{table.1}
\caption{Simulation results of PPCA}
\begin{center}
\begin{tabular}{ccccccccc}
\hline																	
$c=0.4$	&		&		&		&	~	&	$c=2$	&		&		&		\\
\hline																	
$\lambda_1$	&	7.90 	&	$\widehat\lambda_1^{\rm c}$	&	8.37 	&	~	&	$\lambda_1$	&	17.08 	&	$\widehat\lambda_1^{\rm c}$	&	17.13 	\\
$\psi_1$	&	7.92 	&	$\widehat\lambda_1$	&	8.41 	&	~	&	$\psi_1$	&	17.14 	&	$\widehat\lambda_1$	&	17.19 	\\
$\lambda_2$	&	3.75 	&	$\widehat\lambda_2^{\rm c}$	&	3.84 	&	~	&	$\lambda_2$	&	6.27 	&	$\widehat\lambda_2^{\rm c}$	&	5.76 	\\
$\psi_2$	&	3.80 	&	$\widehat\lambda_2$	&	3.89 	&	~	&	$\psi_2$	&	6.43 	&	$\widehat\lambda_2$	&	5.94 	\\
$b$	&	1.21 	&	$\widehat\lambda_{3}$	&	1.17 	&	~	&	$b$	&	2.20 	&	$\widehat\lambda_{3}$	&	2.18 	\\
$a$	&	0.02 	&	$\widehat\lambda_p$	&	0.02 	&	~	&	$a$	&	0.00 	&	$\widehat\lambda_p$	&	0.00 	\\
\hline																																																										
\end{tabular}
\end{center}
\end{table}

\subsection{Ordering-robustness of PPCA}

\begin{figure}[h]
\centering
\includegraphics[height=6cm]{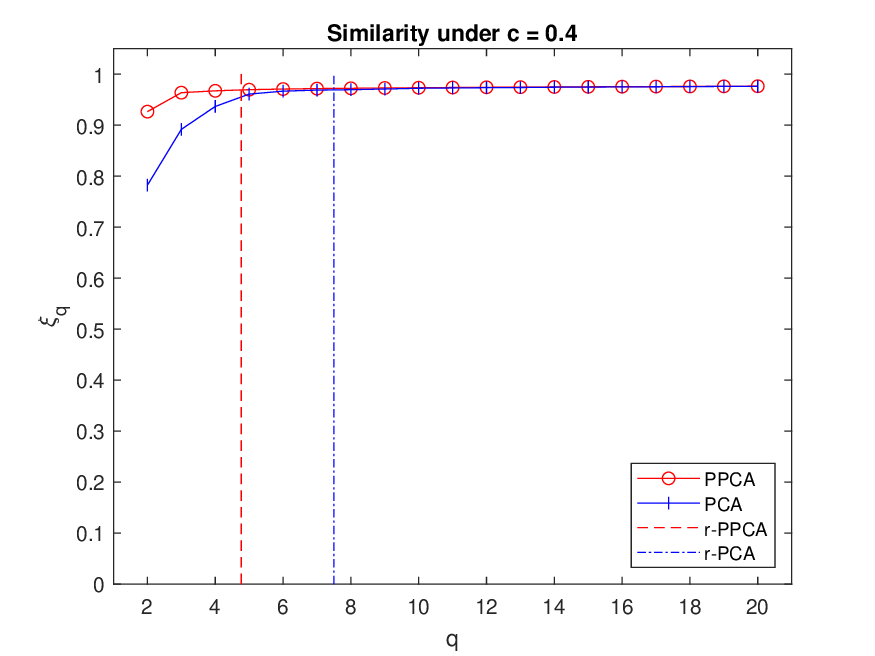}
\includegraphics[height=6cm]{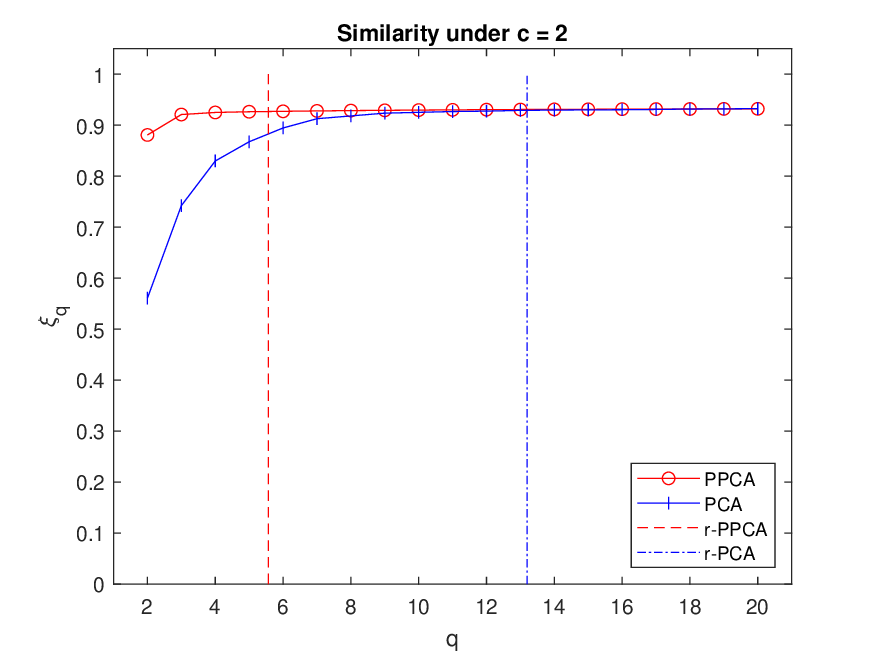}
\caption{The means of similarity $\xi_q$ at $q\in\{2,\ldots, 20\}$ under $c=0.4$ (the left panel) and $c=2$ (the right panel). The vertical lines represent the means of the estimated number of target ranks $\widehat r_{\rm ppca}$ and $\widehat r_{\rm pca}$.}\label{fig.similarity}
\end{figure}

The same simulation setting is considered as in the previous section, except that we set $\nu=2.5$. This indicates that the data set subjects to the presence of heavy-tailed outliers. To measure the influence of outliers on the sample ordering of signal eigenvalues, define the similarity $\xi_q$ as:
\begin{eqnarray}
\xi_q=\frac{1}{r}\sum_{j=1}^{r}s_{qj} \quad {\rm for}\quad q\ge r,
\end{eqnarray}
where $s_{qj}$ denotes the $j^\th$ singular value of $B^\top[\gamma_1,\gamma_2]$, and $B_{p\times q}$ is an orthonormal basis comprising the leading $q$ eigenvectors obtained from PPCA or PCA. A larger value of $\xi_q\in [0,1]$ indicates better performance for $B$ in containing $\Sr$, where $\xi_q=1$ indicates that $\Sr\subseteq\rmspan(B)$. Moreover,
the $\xi_q$ values increases with $q$, and the method with a higher $\xi_q$ curve towards the left-top corner indicates its better performance in preserving the sample ordering of signal eigenvalues. We also estimate the target rank $r_{\rm ppca}$ in (\ref{r.ppca}) and $r_{\rm pca}$ in (\ref{r.pca}) by
\begin{eqnarray}
\widehat r_{\rm ppca}=\sum_{j=1}^{p}\mathcal{I}\{\widehat\lambda_j>b\}\quad{\rm and}\quad \widehat r_{\rm pca}=\sum_{j=1}^{p}\mathcal{I}\{\widetilde\lambda_j>b'\},
\end{eqnarray}
where the upper bounds $b$ and $b'$ are defined in Theorem~\ref{thm.ppca.ssm}. A value of $\widehat r_{\rm ppca}$ (or $\widehat r_{\rm pca}$) deviating from the ground truth $r=2$ to a larger value indicates the tendency of overestimation by PPCA (or PCA) in the presence of outliers. Simulation results (with 200 replicates) are depicted in Figure~\ref{fig.similarity}.

For the case of $c=0.4$, it can be seen that $\widehat r_{\rm PPCA}\approx 5$, and PPCA achieves the $\xi_q$ value around 0.95 at $q=3$. This implies that in the presence of outliers, PPCA tends to have 3 spiked noise eigenvalues, among them one tends to precede the sample ordering of signal eigenvalues. Thus, PPCA cannot recover $\Sr$ until including $\{\widehat\gamma_j\}_{j=1}^3$. PPCA also tends to further include $\{\widehat\gamma_4,\widehat\gamma_5\}$ as signal eigenvectors, although they contribute nearly nothing to the estimation of $\Sr$. On the other hand, it is detected that PCA has a worse performance that PPCA. Specifically, PCA can only achieve the $\xi_q$ value around 0.95 at $q=5$. Together with $\widehat r_{\rm PCA}\approx 8$, it implies that PCA tends to have 6 spiked noise eigenvalues, and three of them tend to precede the sample ordering of signal eigenvalues. As a result, PCA cannot recover $\Sr$ until including $\{\widetilde\gamma_1,\widetilde\gamma_2,\widetilde\gamma_3,\widetilde\gamma_4,\widetilde\gamma_5\}$, and PCA tends to erroneously include $\{\widetilde\gamma_6,\widetilde\gamma_7,\widetilde\gamma_8\}$ in estimating $\Sr$. Comparing the two methods, one can see that PPCA has higher $\xi_q$ values than PCA for $q<6$, and the two curves become close to each other for $q\ge 6$. This supports that, while both PPCA and PCA are able to recover $\Sr$, PPCA produces a more accurate ordering of signal eigenvalues than PCA does. The fact that $\widehat r_{\rm PPCA}<\widehat r_{\rm PCA}$ also supports the robustness of PPCA to the presence of spiked noise eigenvalues, resulting in less bias to the determination of $r$, and hence, better estimation of $\Sr$ than PCA does.

As to the case of $c=2$, a similar conclusion can be made for PPCA and PCA, except that the differences between two $\xi_q$ curves and the gap between $\widehat r_{\rm PPCA}$ and $\widehat r_{\rm PCA}$ become larger than the case of $c=0.4$. In particular, PPCA tends to include $\{\widehat\gamma_j\}_{j=1}^6$ as signal eigenvectors, despite $\{\widehat\gamma_j\}_{j=1}^3$ suffices to recover $\Sr$. On the other hand, PCA cannot recover $\Sr$ until including $\{\widetilde\gamma_j\}_{j=1}^9$, and PCA tends to wrongly include $\{\widetilde\gamma_j\}_{j=9}^{13}$ in the estimation of $\Sr$. This implies that PCA requires to select about 3 more ranks to capture $\Sr$, and tends to produce 7 more spiked noise eigenvalues than PPCA does. Our simulation studies support the robustness of the PPCA eigenvalues to the presence of spiked noise eigenvalues, even in the situation of high-dimensionality.

\section{Discussion}

This paper establishes the high-dimensional asymptotic properties of the PPCA eigenvalues. It is found that PPCA and PCA can have different asymptotic behaviors under the situation of $p\to\infty$ and $\lambda_1<\infty$. This is a distinction to the results of Hung and Huang (2023), where the authors show that PPCA and PCA share the same asymptotic behavior under either the situations of $p<\infty$ and $\lambda_1<\infty$ or $p\to\infty$ and $\lambda_r\to\infty$. Under SSM and in comparison with PCA, while PPCA fails to identify eigenvalues in $[\lambda',\lambda^*]$, PPCA gains more robustness to the presence of spiked noise eigenvalues that can have critical influence to the estimation of $\Sr$. This property supports the usage of PPCA under the high-dimensional setting where outliers are unavoidably encountered.


\section*{References}
\begin{description}
\item
Bai, Z., and Ding, X. (2012). Estimation of spiked eigenvalues in spiked models. Random Matrices: Theory and Applications, 1(02), 1150011.

\item
Bai, Z., and Silverstein, J. W. (2010). {\it Spectral Analysis of Large Dimensional Random Matrices}. 2nd ed, New York: Springer.

\item
Hung, H. and Huang, S. Y. (2023). On the efficiency-loss free ordering-robustness of product-PCA. {\it Submitted}.

\item
Nadakuditi, R. R., and Silverstein, J. W. (2010). Fundamental limit of sample generalized eigenvalue based detection of signals in noise using relatively few signal-bearing and noise-only samples. {\it IEEE Journal of Selected Topics in Signal Processing}, 4(3), 468-480.

\item
Pan, G. (2010). Strong convergence of the empirical distribution of eigenvalues of sample covariance matrices with a perturbation matrix. Journal of Multivariate Analysis, 101(6), 1330-1338.

\item
Yao, J.,  Zheng, S. and Bai, Z. (2015). {\it Large Sample Covariance Matrices and High-Dimensional Data Analysis}. Cambridge University Press.

\item
Yata, K. and Aoshima, M. (2010). Effective PCA for high-dimension, low-sample-size data with singular value decomposition of cross data matrix. {\it Journal of Multivariate Analysis}, 101(9), 2060-2077.
\end{description}

\end{document}